\documentclass[a4paper,11pt]{article} 
\usepackage[english]{babel} 
\usepackage[matrix,arrow,tips,curve]{xy}
\usepackage{amssymb,amsthm}
\usepackage[tbtags]{amsmath}
\usepackage{enumerate}

\newcommand{\pr}[1]{\ensuremath{\mathbb{P}^{#1}}}

\newcommand{\ph}{\ensuremath{\varphi}}
\newcommand{\fx}{\ensuremath{\Sigma_X}}
\newcommand{\fy}{\ensuremath{\Sigma_Y}}
\newcommand{\Z}{\ensuremath{\mathbb{Z}}}
\newcommand{\f}{\ensuremath{\Sigma}}
\newcommand{\Q}{\ensuremath{\mathbb{Q}}}

\newcommand{\C}{\ensuremath{\mathbb{C}}}

\newcommand{\NE}{\operatorname{NE}}
\newcommand{\N}{\ensuremath{\mathcal{N}_1}}
\newcommand{\A}{\ensuremath{\mathcal{A}_1}}
\newcommand{\ld}{.\,.\,}
\newcommand{\w}{\widetilde}

\newcommand{\Hom}{\operatorname{Hom}}
\newcounter{prg}
\newcommand{\prg}[1]{\stepcounter{prg}

\bigskip

{\noindent\bf \theprg.\ {#1}}}

\newtheorem{teo}{Theorem}
\newtheorem*{thm}{Theorem}

\newtheorem{prop}[teo]{Proposition}
\newtheorem{lemma}[teo]{Lemma}

\theoremstyle{definition}
\newtheorem*{defi}{Definition}

\theoremstyle{remark}
\newtheorem*{remark}{Remark}
\newtheorem*{ur}{Useful Remark}

\newtheorem*{claim}{Claim}

\begin{document}
\title{Centrally symmetric generators in toric Fano
  varieties}
\author{Cinzia Casagrande\footnote{Dipartimento di Matematica, 
Universit\`a di Roma ``La Sapienza'',
piazzale Aldo Moro 5, 00185 Roma, e-mail: ccasagra@mat.uniroma1.it}}
\date{22nd October 2002}
{\maketitle}

\bigskip

{\small \noindent {\bf Abstract.}
We give a structure theorem for $n$-dimensional smooth toric Fano
varieties whose associated polytope has ``many'' pairs of centrally
symmetric vertices.}

\bigskip

\noindent{\bf Introduction.}
Smooth toric Fano varieties, together with their equivariant birational
contractions, have been intensively studied in recent
years: see \cite{wat,bat3,bat2,sato,bonavero3,cras,fano,sato2}. In some
cases, the toric case has been a testing ground for results concerning
general smooth Fano varieties, see for instance \cite{bonwisncamp,mukai}.
Despite the strong combinatorial properties of fans of toric Fano
varieties, these varieties are classified only in dimension less or
equal to 4 (see \cite{wat,bat3} for dimension 3, \cite{bat2} and
\cite{sato}, example 4.7, for dimension 4).

Let $X$ be a smooth, complete toric variety of dimension $n$,
$N=\Hom_{\Z}(\C^*,T)$, $N_{\Q}=N\otimes_{\Z}\Q$ and $\fx$ the fan of
$X$ in $N_{\Q}$. We denote by $\langle x_1,\dotsc,x_h\rangle$ the
convex cone spanned in $N_{\Q}$ by $x_1,\dotsc,x_h\in N$.
\begin{defi}
The set of \emph{generators} of $\fx$ is defined in the following way:
\[ G(\fx)=\{x\in N\,|\,\langle x\rangle\in\fx\text{ and 
$x$ is primitive in }\langle x\rangle\cap N\}. \]
\end{defi}
Suppose that $X$ is Fano. Then the  generators $G(\fx)$ are the vertices of a
simplicial convex polytope $\mathcal{P}_X\subset N_{\Q}$, and  every cone in
$\fx$ is the cone over a face of $\mathcal{P}_X$. Hence the set of
generators $G(\fx)$
determines the fan $\fx$ (and $X$). 
This means that $X$ is the only Fano
variety in the set
\begin{align*} 
\mathcal{B}_X=\{\,& X' \text{ toric, smooth, complete }|\,\exists\, \ph\colon
X'\dasharrow X\text{ equivariant, birational, isomorphism} \\
&\text{in codimension 1}\}\\
=\{\,& X' \text{ toric, smooth, complete }|\ G(\f_{X'})=G(\fx)\,\}
\end{align*}
and any such $X'$ determines $X$. It is an interesting question
whether there exists some $X'$ in $\mathcal{B}_X$ with a particularly
simple structure. More precisely,
we think of a toric bundle structure (see section 2), 
with fiber and basis some lower
dimensional toric Fano variety. This would give an inductive way to
classify these varieties.
Our result is the following:
\begin{thm}[theorem \ref{cs} and
proposition \ref{dipiu}]
Let $X$ be a smooth toric Fano variety.
Suppose that $G(\fx)$ contains $r>0$ pairs of centrally symmetric
generators $\{x,-x\}$, and
consider the linear subspace $H$ of
$N_{\Q}$ spanned by all such pairs. 
\begin{enumerate}[$\bullet$]
\item If $\dim H<r$, then $X$ is a toric bundle over a lower
  dimensional toric Fano variety, with 
fiber a product of Del Pezzo  varieties (see
definition on page \pageref{pseudo}). 
\item 
If $H$ contains more than $2\dim H$ generators, 
then there exists $X'$ in $\mathcal{B}_X$ which
is a toric bundle over a lower dimensional toric Fano variety, with
fiber a product of Del Pezzo and pseudo Del Pezzo varieties (see
definitions on page \pageref{pseudo}).
\item If $\dim H=n-1$, then there exists $X'$ in
$\mathcal{B}_X$ which 
is a toric bundle over $\pr{1}$, with
fiber a product of Del Pezzo varieties, pseudo Del Pezzo varieties and
a power of $\pr{1}$.
\item If  $\dim H=n$, then $X$ itself is
 a product of Del Pezzo varieties, pseudo Del Pezzo varieties and a
 power of $\pr{1}$.
\end{enumerate}
\end{thm}
All these results give a structure theorem for some classes of toric
Fano varieties, in any dimension. The last case ($\dim H=n$)
generalizes a previous result by G.~Ewald~\cite{ewald},
who obtaines the same statement under the stronger hypothesis that
 there exists an $n$-dimensional cone $\sigma\in\fx$
such that $-\sigma\in\fx$.

The structure of the paper is as follows: in section 1 we recall the 
language of primitive relations and some results of toric Mori theory; 
in section 2 we recall properties of toric bundles; in section 3 we 
introduce Del Pezzo and pseudo Del Pezzo varieties and state and prove
our results.

\prg{Primitive relations and 1-cycles.} We collect here some basic
notions and results. We refer to \cite{oda,fulton} for properties of toric
varieties, to \cite{bat2} for
primitive collections and primitive relations, to \cite{reid} and to
the survey paper \cite{torimori} for
toric Mori theory.

Let $X$ be a smooth, complete toric variety.
\begin{defi}
A set $P\subseteq G(\fx)$ is a \emph{primitive collection} if $\langle
P\rangle\not\in\fx$ and  $\langle
P\smallsetminus\{x\}\rangle\in\fx$ for all $x\in P$.

If $P=\{x_1,\dotsc,x_h\}$ is a primitive collection, consider the
minimal cone $\sigma\in\fx$ containing the point $x_1+\cdots+x_h$. If
$\sigma=\langle y_1, \dotsc,y_k\rangle$, we get a relation
\[ x_1+\cdots+x_h-(a_1y_1+\cdots+ a_ky_k) =0\]
with $a_i\in\Z_{>0}$. This is the \emph{primitive relation} $r(P)$ associated
to $P$, and the \emph{degree} of $P$ is $\deg P=h-\sum_ia_i$.
\end{defi}
Let $\A(X)$ be the group of algebraic 1-cycles on $X$ modulo numerical 
equivalence, $\N(X)=\A(X)\otimes_{\Z}\Q$, and $\NE(X)\subset \N(X)$ the
convex cone generated by classes of effective curves; we recall that
$\NE(X)$ is closed and polyhedral.

For every $r\in\{0,\dotsc,n\}$ there is a bijection between
the set of $r$-dimensional cones in $\fx$ and the 
set of codimension $r$ orbits of $\,T$ in $X$. For $r=1$, this gives a
bijection between $G(\fx)$   and the 
set of codimension 1 orbits of $\,T$. We denote by $V(\sigma)$ 
(respectively $V(x)$)
 the closure of the orbit corresponding to $\sigma\in\fx$ 
(respectively $x\in
 G(\fx)$); this is an irreducible
invariant subvariety.

We recall that there is an exact sequence:
\[
0\longrightarrow \A(X)\longrightarrow
\mathbb{Z}^{G(\fx)}\longrightarrow N \longrightarrow 0 \]
where the map $\A(X)\rightarrow \mathbb{Z}^{G(\fx)}$ is given by 
$\gamma\mapsto \{\gamma\cdot V(x)\} _{x\in G(\fx)}$. Hence the group 
$\A(X)$ is 
canonically isomorphic to the lattice of
integral relations among the elements of $G(\fx)$:
a relation $\sum_x a_x x=0$ corresponds to a numerical
class $\gamma$ that has
intersection $a_x$ with each divisor 
$V(x)$. Moreover, we have $(-K_X)\cdot\gamma=
\sum_xa_x$. In particular, every primitive relation $r(P)$ defines a
numerical class in $\A(X)$, and $\deg P=(-K_X)\cdot r(P)$.
In the sequel, we will identify numerical classes with their
associated relations.
Remark that from the exact sequence above we deduce
$\rho_X=\#G(\fx)-n$, where $\rho_X$ is the Picard number of $X$.

\begin{lemma}[\cite{contr}, lemma 1.4]
\label{lemma}
Let $\gamma\in\A(X)$ given by the relation
\[ a_1x_1+\cdots+a_hx_h-(b_1y_1+\cdots+b_ky_k)=0
\]
with $a_i,b_j\in\Z_{>0}$ for each $i,j$. If $\langle
y_1,\dotsc,y_k\rangle\in\fx$, then $\gamma\in\NE(X)$.
\end{lemma}
In particular, every primitive relation is an effective class:
therefore, \emph{if $X$ is Fano, every primitive collection has strictly
positive degree}.
\begin{teo}
\label{unoextr}
Let $X$ be a toric Fano variety and
let $\gamma\in\NE(X)\cap\A(X)$ 
be such that $(-K_X)\cdot\gamma=1$. Then $\gamma$
is extremal in $\NE(X)$, the relation associated to $\gamma$ is a
primitive relation
\[x_1+\cdots+x_h-(a_1y_1+\cdots+a_ky_k)=0,
\]
and for every cone $\nu=\langle z_1,\dotsc,z_t\rangle$ 
such that $\{z_1,\ld,z_t\}\cap\{ x_1,\ld,x_h,y_1,\ld,y_k\}=\emptyset$
and $\langle y_1,\ld,y_k\rangle+\nu \in\fx$,
we have $\langle
x_1,\ld,\check{x}_i,\ld,x_h,y_1,\ld,y_k\rangle+\nu
\in\fx$ for all $i=1,\dotsc,h$.
\end{teo}
\begin{proof}
See \cite{contr}, corollary 4.4 and \cite{reid}, theorem 2.4.
\end{proof}
\begin{lemma}[\cite{fano}, lemma 2.2]
\label{idea}
Let $X$ be a toric Fano variety.
Suppose that $\{x,-x\}$ is a
primitive collection in $\fx$. Then for any other primitive collection $P$
containing $x$, the associated primitive relation is  
\[r(P):\quad x+y_1+\cdots+y_h-(z_1+\cdots+z_h)=0; \] 
moreover, $P'=\{-x,z_1,\dotsc,z_h\}$ is also a primitive collection,
with relation  
\[r(P'):\quad  (-x)+z_1+\cdots+z_h-(y_1+\cdots+y_h)=0. \]
Both $P$ and $P'$ have degree 1, their primitive relations 
are extremal, and $2h\leq n$.
\end{lemma}

\vspace{-6pt}

\prg{Toric bundles.}
\label{bundle}
We recall the following standard result of toric
geometry (see~\cite{ewald2}, theorem 6.7 in chapter VI,
or~\cite{fulton}, section 2.4, exercise on fiber bundles):
\begin{teo}
\label{soleil}
Let $X$ be a smooth, complete toric variety. Consider a linear
subspace $H\subset N_{\Q}$ of dimension $r$ and suppose that for every
$n$-dimensional cone $\sigma\in\fx$, we have $\sigma=\eta+\tau$,  with 
$\eta\in\fx$, $\eta\subset H$, $\dim\eta=r$ and $\tau\cap H=\{0\}$. Then: 

\smallskip

\noindent $\bullet\ $ the set $\f_F=\{\sigma\in\fx\,|\,\sigma\subset
H\}$ is the
fan of a smooth, complete, $r$-dimensional toric variety $F$;

\smallskip

\noindent $\bullet\ $ consider a complementar $H^{\perp}$ 
of $H$ in $N_{\Q}$ such
that $N=(N\cap H)\oplus(N\cap H^{\perp})$, and let $\pi\colon
N_{\Q}\rightarrow H^{\perp}$ be the projection. Then the set $\f_Z=\{\pi
( \sigma)\,|\,\sigma\in\fx\}$ is the fan of a complete, smooth,
$(n-r)$-dimensional toric variety $Z$;

\smallskip

\noindent $\bullet\ $ the projection $\pi$ induces an
equivariant morphism $\overline{\pi}\colon X\rightarrow Z$ such that, for every
affine invariant open subset $U\subset Z$, $\overline{\pi}^{\,-1}(U)\simeq 
F\times U$ as toric varieties over $U$. 

\smallskip

Hence, $X$ is a toric bundle over $Z$ with fiber $F$.
\end{teo}
Under this description, for every 
$\eta\in\fx$ such that $\eta\subset H$ and $\dim\eta=r$, 
the subvariety
$V(\eta)$ is an invariant section 
of the bundle; for every $\tau\in\fx$
such that $\tau\cap H=\{0\}$ and $\dim\tau=n-r$, 
the subvariety $V(\tau)$ is an invariant fiber.

Suppose that $X$ is a toric bundle over $Z$ with fiber $F$, as
described in theorem~\ref{soleil}. 
Then  $G(\f_F)=G(\fx)\cap H$ and
the primitive collections of $\fx$ contained in $H$ are exactly the
primitive collections of $\f_F$. These primitive collections have the
same primitive relations in $\fx$ and in $\f_F$.

Now set $\mathcal{G}=G(\fx)\smallsetminus G(\f_F)$. The
projection $\pi\colon N_{\Q}\rightarrow H^{\perp}$ 
induces a bijection between $\mathcal{G}$ and
$G(\f_Z)$. 
Under this identification, the primitive collections of
$\fx$ not contained in $H$ are exactly the primitive collections of
$\f_Z$. We denote by $\overline{x}$ the image of an element $x\in N$
under $\pi$. A primitive relation in $\f_Z$ of the form 
\[ \overline{z}_1+\cdots+\overline{z}_h-(a_1 \overline{y}_1+\cdots+a_k
\overline{y}_k)=0\]
will lift in $\fx$ to a primitive relation of the form:
\[ z_1+\cdots+z_h-(a_1 y_1+\cdots+a_ky_k+b_1x_1+\cdots+b_lx_l)=0,\]
with $l\geq 0$ and $x_i\in G(\f_F)$ for all $i=1,\dotsc,l$. 

The bundle is trivial, namely $X\simeq F\times Z$, 
if and only if all primitive relations in
$\f_Z$ remain unchanged when lifted in $\fx$.

\prg{Del Pezzo varieties and centrally symmetric generators.} 
Let $x_1,\dotsc,x_n$ be a basis of the lattice $N$ and suppose that
$n=2r$ is even. Consider the following elements of $N$:
\[ x_0=-x_1-\cdots-x_n,\qquad y_i=-x_i\quad\text{ for }i=0,\dotsc,n.\]

\noindent\emph{Notation:} $I^s$ denotes a set of indexes of
cardinality $s$.
\begin{defi}
\label{pseudo}
The $n$-dimensional \emph{Del Pezzo variety $V^n$} is the toric Fano variety
described by the fan $\f_{V^n}\subset N_{\Q}$ such that
$G(\f_{V^n})=\{x_0,y_0,\ld,x_n,y_n\}$. 

The $n$-dimensional \emph{pseudo Del Pezzo variety $\widetilde{V}^n$} is the
toric Fano variety described by the fan $\f_{\widetilde{V}^n}\subset
N_{\Q}$ such that
$G(\f_{\widetilde{V}^n})=\{x_0,x_1,y_1,\ld,x_n,y_n\}$.
\end{defi}
We have $\rho_{V^n}=n+2$ and $\rho_{\w{V}^n}=n+1$.
In dimension 2, these varieties are 
the Del Pezzo surfaces $S_3$ and $S_2$, obtained
from $\pr{2}$ blowing-up respectively three or two
fixed points. If the dimension is greater
than two, 
the varieties $V^n$ and $\widetilde{V}^n$ do not admit any smooth equivariant
blow-down;
 moreover, the blow-up along any invariant subvariety is
never Fano. This behaviour is rather peculiar: in dimension 4, 
H.~Sato~(\cite{sato}, theorem 8.1)
has shown that any toric Fano 4-fold different from $V^4$ and
$\widetilde{V}^4$ can be obtained from $\pr{4}$ by a sequence of smooth
equivariant blow-ups and blow-downs in such a way that all the
intermediate varieties are Fano.

The variety $V^n$ can be geometrically described as follows. Consider
the Cremona map $\pr{n}\dasharrow\pr{n}$. It is well known that 
this map factorizes
as a sequence of smooth blow-ups and blow-downs in the following way:
first blow-up $n+1$ general points $\,p_1,\dots,p_{n+1}\,$ in $\pr{n}$. 
Then, for all
$m=2,\dotsc,n-2$ (in this order), flip the strict transforms of all
the $\pr{m-1}$ spanned by $m$ of the points $p_i$.
Finally, blow-down the strict transforms of all
the $\pr{n-1}$ spanned by $n$ of the points $p_i$.
It is easy to see that this is a toric
factorization, namely, all varieties involved are toric and the maps
are equivariant. Now, $V^n$ is the intermediate variety that
you get if you perform the flips just till the dimension $m=r-1$,
where $r=n/2$. 

Del Pezzo varieties were introduced by V.~E.~Voskresenski{\u{\i}}
and A.~A.~Klyachko \cite{VK}, where they show (theorem 6)
that \emph{every toric Fano
variety whose fan is centrally symmetric is a product of Del Pezzo
varieties and a power of $\pr{1}$.} This result has been generalized by
G.~Ewald~\cite{ewald}, who
introduces pseudo Del Pezzo varieties and shows that 
\emph{every toric Fano
variety whose fan contains two centrally
symmetric maximal cones, is a product of Del Pezzo
varieties, pseudo Del Pezzo varieties and a power of $\pr{1}$.} 
We will see in proposition~\ref{dipiu} that it is actually enough to ask
that $\fx$ has $n$ linearly independent pairs of symmetric generators.

Explicitly, the fans of ${V^n}$ and ${\w{V}^n}$ are:
\begin{align*}
\f_{V^n}=\{&\langle x_i,y_j\rangle_{i\in I^r,\,j\in J^r}\,|\,I^r,J^r\subset
  \{0,\dotsc,n\}\text{ disjoint}\}; \\
\f_{\w{V}^n}= \{&\langle x_0,x_i,y_j\rangle_{i\in
  I^{r-1},\,j\in J^r},  
\langle x_i,y_j\rangle_{i\in\w{I}^{r+s},\,j\in \w{J}^{r-s}}\,| \,
I^{r-1},J^r\subset
  \{1,\dotsc,n\}\text{ disjoint},\\
& s\in\{0,\dotsc,r\}\text{ and }
 \w{I}^{r+s},\w{J}^{r-s}\text{ a partition
  of }\{1,\dotsc,n\}\}.
\end{align*}
The primitive relations of $\f_{V^n}$ are $x_i+y_i=0$ for $i=0,\dotsc,n$ and
\[ \sum_{i\in I^{r+1}} x_i-(\sum_{j\in J^r}y_j)=0,\qquad\sum_{i\in
  I^{r+1}} y_i-(\sum_{j\in J^r}x_j)=0
\] for every partition $I^{r+1},J^r$ of $\{0,\dotsc,n\}$.

The primitive relations of $\f_{\w{V}^n}$ are $x_i+y_i=0$ for $i=1,\dotsc,n$ and
\[ x_0+\sum_{i\in I^{r}} x_i-(\sum_{j\in J^r}y_j)=0,\qquad\sum_{j\in
  J^{r+1}} y_j-(x_0+\sum_{i\in I^{r-1}}x_i)=0
\] for all partitions $I^{r},J^r$ and $I^{r-1},J^{r+1}$ of
$\{1,\dotsc,n\}$. 

What makes Del Pezzo 
and pseudo Del pezzo
varieties special among toric Fano varieties, is that they have lots of
numerical classes of curves with anticanonical degree 1.
\begin{teo}
\label{cs}
Let $X$ be a toric Fano variety and let $H$ be the linear subspace of
$N_{\Q}$ spanned by all pairs of centrally symmetric generators in
$G(\fx)$. 
Then there exists a decomposition
\[ H=H_1\oplus\cdots\oplus H_a\oplus\widetilde{H}_1\oplus\cdots\oplus
\widetilde{H}_b\oplus K \]
with $a\geq 0$, $h_i=\dim H_i>0$, $b\geq 0$,
$\w{h}_i=\dim\widetilde{H}_i>0$ and $\dim K\geq 0$, such that:
\begin{enumerate}[(i)]
\item 
$G(\fx)\cap H= \bigcup_{i=1}^a 
\left(G(\fx)\cap H_i \right)\cup \bigcup_{j=1}^b
(G(\fx)\cap \widetilde{H}_j )\cup \left(G(\fx)\cap K\right)$;
\item $\#(G(\fx)\cap H_i)=2(h_i+1)$, $\#(G(\fx)\cap \w{H}_j)=2\w{h}_j+1$
  and $\#(G(\fx)\cap K)=2\dim K$ for all $i=1,\dotsc,a$ and $j=1,\dotsc,b$;
\item there exists a commutative diagram 
\[
\xymatrix{
{X} \ar[d] \ar@{-->}[r]^{\ph}  &  {X'} \ar[d]\ar[dr] & \\
Z \ar@{-->}[r]^{\psi} & {Z'} \ar[r] & W }
\]
where all varieties are smooth, projective and toric, 
$Z$ and $W$ are Fano,  $\ph$ and $\psi$ are
equivariant birational maps which are
 isomorphisms in codimension 1, and the decomposition of $H$ induces
toric bundle structures:
\begin{itemize}
\item  $X\rightarrow Z$  and  $X'\rightarrow Z'$ with fiber
$V^{h_1}\times\cdots\times V^{h_a}$, 
\item
$Z'\rightarrow W$ with fiber
  $\w{V}^{\w{h}_1}\times\cdots \times\w{V}^{\w{h}_b}$, 
\item
$X'\rightarrow W$ with fiber
$V^{h_1}\times\cdots\times V^{h_a}\times\w{V}^{\w{h}_1}\times\cdots 
\times\w{V}^{\w{h}_b}$.
\end{itemize}
\end{enumerate}
\end{teo}

\noindent\emph{Remarks:} 
\begin{enumerate}[$\bullet$] 
\item $X'\in\mathcal{B}_X$.
\item
in particular, we will see in the proof 
(step 7) that $\ph$ is an isomorphism over
some invariant section of $X'\rightarrow W$, hence
there exists an invariant subvariety $V\simeq W$ of $X$:
$V$ is Fano, $\dim V=
n-\sum_i h_i-\sum_j\w{h}_j$ and $\rho_{V}=\rho_X-\sum_i (h_i+2)
-\sum_j (\w{h}_j+1)$ (where $\rho_V$ and $\rho_X$ are the Picard
number of $V$ and $X$ respectively).
\item when $H_i$ or $\w{H}_i$ has dimension 2, we
recover \cite{fano}, 2.4.3 and 2.4.2a respectively. 
\item for all $i=1,\dotsc,a$,
$G(\fx)\cap H_i$ consists of $h_i+1$ pairs of symmetric generators
$\{x_l,y_l\}$, $h_i$ by $h_i$ linearly independent, such that
$x_1+\cdots+x_{h_i+1}=0$. 
\item for all $j=1,\dotsc,b$, $G(\fx)\cap \w{H}_j$ 
consists of $\w{h}_j$ linearly independent 
pairs of symmetric generators
$\{x_l,y_l\}$, plus another generator
$x_0=x_1+\cdots+x_{\w{h}_j}$.
\item $G(\fx)\cap K$ consists of $\dim K$ linearly
independent pairs of symmetric generators. 
\end{enumerate}

For the proof of theorem \ref{cs}, we need the following lemma:
\begin{lemma}
\label{dopo}
Suppose that $\fx$ contains $m$ pairs of centrally symmetric generators
$x_i+y_i=0$ for $i=1,\dotsc,m$. Let $\tau=\langle
z_1,\dotsc,z_s\rangle\in\fx$ such that
$\{z_1,\ld,z_s,x_1,\ld,x_m\}$ are linearly independent. Then
there exists a partition $I,J$ of $\{1,\ld,m\}$ such that
$\tau+\langle
x_i,y_j \rangle_{i\in I,\,j\in J}\in\fx$.
In particular, $\{z_1,\ld,z_s,x_1,\ld,x_m\}$ is a part of a basis of
the lattice $N$.
\end{lemma}
\begin{proof}[Proof of lemma \ref{dopo}]
Suppose that $\tau+\langle y_1\rangle\not\in\fx$. Then, up to
 renumbering $\{z_1,\ld,z_s\}$,
 there is a
primitive collection $\{z_j,\dotsc,z_s,y_1\}$ with $1\leq j\leq s$, and
by lemma~\ref{idea} we get primitive relations:
\[ z_j+\cdots+z_s+y_1-(u_1+\cdots+u_{s-j+1})=0,\quad u_1+\cdots+u_{s-j+1}+x_1-
(z_j+\cdots+z_s)=0.\]
Now, theorem~\ref{unoextr}
implies that either $\tau+\langle x_1\rangle\in\fx$, or
$\{z_1,\dotsc,z_s\}\supseteq \{u_1,\dotsc,u_{s-j+1}\}$. This last case is
impossible, because it would give a relation of linear dependence
among $\{z_1,\ld,z_s,y_1\}$, against the hypotheses.
Therefore $\tau+\langle x_1\rangle\in\fx$, and by recurrence we have
the statement. 
\end{proof}
\begin{ur}
Suppose that we have a relation among generators:
\[ a_1x_1+\cdots+a_hx_h-(b_1y_1+\cdots+b_ky_k)=0
\]
with $a_i,b_j\in\Z_{>0}$ for each $i,j$. 
Then the two cones $\langle x_1,\dotsc,x_h\rangle$ and $\langle 
y_1,\dotsc,y_k\rangle$ can not be both in $\fx$, because they
intersect in their interior, while two cones in the fan always
intersect along a face.
\end{ur}
\begin{proof}[Proof of theorem \ref{cs}]
We fix an $m$-dimensional linear subspace $H'$ of $H$ spanned by $m$
pairs of centrally symmetric generators $x_i+y_i=0$,
with  $x_1,\dotsc,x_m$ linearly independent.

\smallskip

\noindent\emph{Step 1:}
every generator of $\fx$ contained in $H'$ and different from
$x_1,y_1,\dotsc,x_m,y_m$
 has the form $\sum_{i=1}^m
\varepsilon_i x_i$ with $\varepsilon_i\in\{1,0,-1\}$ for all $i$.
\begin{proof}[Proof of step 1]
Let $x_0 \in G(\fx)$ be a generator in $H'$ different from
$\{x_1,y_1,\dotsc,x_m,y_m\}$. Up to reorder $\{x_1,\dotsc,x_m\}$, we can
suppose that $x_0=\sum_{i=1}^s \varepsilon_i x_i$, with
$2\leq s\leq m$ and
$\varepsilon_i\neq 0$ for all $i=1,\dotsc,s$. 
Then for all $j\in\{1,\ld,s\}$  
the generators
$\{x_0,x_1,\ld,\check{x}_j,\ld,x_s\}$
are linearly independent, and 
by lemma~\ref{dopo} they are a part of a basis of the
lattice, hence $|\varepsilon_j|=1$.
\end{proof}
We  assume from now on that
 $x_0=y_1+\cdots+y_m \in G(\fx)$.

\smallskip

\noindent \emph{Step 2: } $m=2r$ is
even and $\fx$ contains all the cones 
\begin{equation}
\tag{$\star$}
\label{coni}
\langle x_i,y_j\rangle_{i\in I^r,\,j\in J^r},\qquad 
  \langle x_0,x_i,y_j\rangle_{i\in \widetilde{I}^{r-1},\,j\in 
\widetilde{J}^r} 
\end{equation}
and the primitive extremal relations
\[
x_0+\sum_{i\in I^r}x_i-(\sum_{j\in J^r}y_j)=0,\qquad \sum_{J^{r+1}}y_j
-(x_0+\sum_{i\in I^{r-1}} x_i)=0 \]
for all partitions $I^r,J^r$ and $I^{r-1},J^{r+1}$ of $\{1,\dotsc,m\}$ 
and for all choices of $\widetilde{I}^{r-1},
\widetilde{J}^r\subset\{1,\dotsc,m\}$ disjoint.

Moreover, for every $\tau\in\fx$ such that $\tau\cap H'=\{0\}$ and 
$\tau+\eta\in\fx$ for some $\eta$ in \eqref{coni}, then
$\tau+\eta'\in\fx$ for all $\eta'$ in \eqref{coni}.
\begin{proof}[Proof of step 2]
By lemma~\ref{dopo}, since $\{x_0,x_1,\ld,x_{m-1}\}$ are linearly
independent and $\langle x_0\rangle\in\fx$, 
up to reorder $\{x_1,\dotsc,x_{m-1}\}$
we have  $\langle
x_0,x_1,\ld,x_{r-1},y_r,\ld,y_{m-1}\rangle\in\fx$ with $1\leq r\leq m$. 
Applying again lemma~\ref{dopo} to the cone
$\sigma=\langle x_1,\ld,x_{r-1},y_r,\ld,y_{m-1}\rangle$, we get
that at least one of the two cones $\sigma+\langle y_m\rangle$,
$\sigma+\langle x_m\rangle$ is in $\fx$. 

Consider the relation $y_r+\cdots+y_m-(x_0+x_1+\cdots+x_{r-1})=0$. Since 
$\langle x_0,x_1,\ld,x_{r-1}\rangle\in\fx$, by lemma~\ref{lemma}
this relation corresponds to an effective class, and since $X$ is
Fano, its degree must be positive: this gives $2r\leq m$.
Moreover, it can not be $\langle y_r,\dotsc,y_m\rangle\in\fx$
(see the Useful Remark), hence 
$\sigma+\langle x_m\rangle=\langle x_0,x_1,\ld,x_{r-1},y_r,\ld,y_{m-1},x_m
\rangle\in\fx$. This implies that also the relation
$x_0+x_1+\cdots+x_{r-1}+x_m-(y_r+\cdots+y_{m-1})=0$ is effective, hence
its degree is positive, and we get $2r\geq m$.
Therefore $2r=m$. In particular, the relation
$x_0+x_1+\cdots+x_{r-1}+x_m-(y_r+\cdots+y_{m-1})=0$ has degree 1, so 
by theorem~\ref{unoextr}
$\langle x_1,\ld,x_{r-1},x_m,y_r,\ld,y_{m-1}\rangle\in\fx$. 

Now, we know that at least one cone in each of the
two types in (\ref{coni}) is in $\fx$. 
Let $I^r$, $J^r$ be such that the cone 
$\langle x_i,y_j\rangle_{i\in
  I^r,\,j\in J^r}$ is in $\fx$. 
 Then the relation
\[ x_0+\sum_{i\in I^r} x_i -(\sum_{j\in J^r} y_j)=0 \]
is effective by lemma~\ref{lemma}, thus
theorem~\ref{unoextr} implies that the cones 
\[\langle x_0,x_i,y_j\rangle_{i\in I^r\smallsetminus\{p\},\,
j\in  J^r}\] are in $\fx$ for all $p\in I^r$.

In an analogue way, suppose that  for some $\widetilde{I}^{r-1}$, 
$\widetilde{J}^r$ the cone 
$\langle x_0,x_i,y_j\rangle_{i\in
  \widetilde{I}^{r-1},\,j\in \widetilde{J}^r}$ is in $\fx$.
Let $\{1,\dotsc,m\}=\widetilde{I}^{r-1}\cup \widetilde{J}^r\cup\{l\}$. 
 Then the relations
\[\sum_{j\in \widetilde{J}^r\cup\{l\}} 
y_j -(x_0+\sum_{i\in \widetilde{I}^{r-1}} x_i)=0,\qquad
x_0+ \sum_{i\in\widetilde{I}^{r-1}\cup\{l\}}x_i-(\sum_
{ j\in\widetilde{J}^r} y_j)=0
\]
are effective by lemma~\ref{lemma}, and 
again,
theorem~\ref{unoextr} implies that the cones 
\[\langle x_0,x_i,y_j\rangle_{i\in \widetilde{I}^{r-1},\,
j\in(\widetilde{J}^r\smallsetminus\{p\})\cup\{l\}},\quad
\langle x_i,y_j\rangle_{i\in \widetilde{I}^{r-1}\cup\{l\},\,
j\in\widetilde{J}^r}
\]
 are in $\fx$ for all $p\in \widetilde{J}^r$.
Hence, proceeding in this way, we get all the desired cones 
and extremal relations in $\fx$.

Finally, if $\tau$ is such that $\tau\cap H'=\{0\}$ and 
$\tau+\eta\in\fx$ for some $\eta$ in \eqref{coni}, then it is easy to
see, using the extremal relations and reasonning as above, that
$\tau+\eta'\in\fx$ for all $\eta'$ in \eqref{coni}.
\end{proof}
\noindent \emph{Step 3:} set $y_0=-x_0$. Either $H'\cap
G(\fx)=\{x_0,x_1,y_1,\ld,x_m,y_m\}$, or 
$H'\cap G(\fx)=\{x_0,y_0,x_1,y_1,\ld,x_m,y_m\}$. 
\begin{proof}[Proof of step 3]
Suppose that $z\in H'\cap G(\fx)$,
$z\not\in\{x_0,x_1,y_1,\ld,x_m,y_m\}$: we want
to show that $z=y_0=-x_0$.
We know by step 1
that $z=\sum_{i=1}^m\lambda_i y_i$ with $\lambda_i\in\{1,0,-1\}$ for all
$i$:
\[ z=\sum_{\lambda_i=1}y_i+\sum_{\lambda_i=-1}x_i=x_0+
\sum_{\lambda_i=0}x_i+\sum_{\lambda_i=-1}2x_i. \]
Since $z$ can not be in the interior of some cone of
$\fx$, 
we know by step 2 that  either $\#\{i\,|\,\lambda_i=1\}\geq r+1$ or 
$\#\{i\,|\,\lambda_i=-1\}\geq r+1$,
and that $\#\{i\,|\,\lambda_i=0\text{ or
  }\lambda_i=-1\}\geq r$. Hence at least $r+1$ of the $\lambda_i$
are -1.
Thus we have
\[ x_0+z=\sum_{i=1}^m (\lambda_i+1)y_i=\sum_{\lambda_i=1} 2y_i
+\sum_{\lambda_i=0}y_i \]
and since there are at most $r-1$ $y_i$ appearing on the right hand side,
we know by step 2 that they generate a cone in $\fx$. Hence
$\{x_0,z\}$ is a
primitive collection (see the Useful Remark), and since $X$ is Fano, 
up to reordering $\{y_1,\dotsc,y_m\}$
there are only two possible primitive relations: either $x_0+z=0$ 
(so $z=y_0$), 
or $x_0+z-y_1=0$ (so $z=x_2+\cdots+x_m$). Let's show that this last case
is impossible.  
If $x_0+z-y_1=0$, this relation has degree 1; thus by theorem~\ref{unoextr}
$\langle y_1,\ld,y_r\rangle\in\fx$ implies $\langle
z,y_1,\ld,y_r\rangle\in\fx$. But we also have $z+y_2+\cdots
+y_r=x_{r+1}+\cdots+x_m$ and $\langle x_{r+1},\ld,x_m\rangle\in\fx$, a
contradiction (by the Useful Remark). 
\end{proof}

\noindent\emph{Step 4:} for every generator $w\in G(\fx)$ such that $w\not\in
H'$, we have $\langle w \rangle +\eta\in\fx$ for all $\eta$ in
\eqref{coni}.
\begin{proof}[Proof of step 4]
  Suppose first that $\langle w,x_0\rangle\in\fx$. Since  
$\{w,x_0,x_1,\ld,x_{m-1}\}$ are linearly independent,
we know by lemma~\ref{dopo} that there exists a partition $I,J$ of
$\{1,\ld,m-1\}$ such that $\langle w,x_0,x_i,y_j\rangle_{i
  \in I,\,j\in J} \in\fx$. This cone can not contain a primitive
collection, so it must be $\#I\leq r-1$, $\#J\leq r$. Since
$\#I+\#J=2r-1$, it must be $\#I=r-1$, $\#J=r$. Hence the statement
follows from step 2.

If $\langle w,x_0\rangle\not\in\fx$, it means that $\{w,x_0\}$ is 
a primitive collection. Since $w$ is not in $H'$, $w\neq -x_0$; hence the only
possible primitive relation is $w+x_0-u=0$, which has degree 1 and is
extremal by theorem~\ref{unoextr}. Now, $u\not\in H'$ and
$\langle u,x_0\rangle\in\fx$: so by what preceeds, $\langle
u\rangle+\eta\in\fx$ for all $\eta$  in \eqref{coni}. 
Choose $\eta=\langle x_i, y_j\rangle_{i\in I^r,\,j\in J^r}$ for some
$I^r$, $J^r$. Since the
relation  $w+x_0-u=0$ is extremal, theorem~\ref{unoextr} implies
that also $\langle w\rangle+\eta$ is in $\fx$. By step 2, the same is
true for all $\eta'$ in \eqref{coni}, thus we get $\langle
w,x_0\rangle\in\fx$, a contradiction.
\end{proof}

\noindent \emph{Step 5:} let $W$ be the invariant subvariety 
 associated to one of the cones in (\ref{coni}). 
Then $W$ is Fano. 
\begin{proof}[Proof of step 5]
We remark first of all that the linear span of the cones in
\eqref{coni} is always $H'$, and by step 2 the invariant subvarieties
associated to these cones are all isomorphic. Let's fix a cone
$\eta$ in \eqref{coni} and let
$W=V(\eta)$. 

Consider an invariant curve $C$ in $W$: in $\fx$ 
the curve $C$ is associated to the
$(n-1)$-dimensional cone $\langle
w_3,\dotsc,w_s\rangle+\eta=(\langle
w_1,w_3,\ld,w_s\rangle+\eta) \cap(\langle w_2,w_3,\ld,w_s\rangle+\eta)$.
Remark that by step 2, $\langle w_1,w_3,\ld,w_s\rangle+\eta'\in\fx$ and 
$\langle w_2,w_3,\ld,w_s\rangle+\eta'\in\fx$ for all $\eta'$ in \eqref{coni}.
The numerical class of $C$ in $W$ is given by
the relation
\[ \sum_{i=1}^s a_i \overline{w}_i=0, \]
where $a_i\in\mathbb{Z}$ and $\overline{w}$ denotes the image of $w\in
N$ under the projection $N\rightarrow N/H'$. Hence in $W$ the
anticanonical degree of $C$ is $\sum_i a_i$.
In the lattice $N$, this relation lifts as
\[ \sum_{i=1}^s a_iw_i +\sum_{j=1}^m \lambda_jx_j=0, \]
where $\lambda_j\in\mathbb{Z}$. 
Choose a partition $I^r$, $J^r$ of $\{1,\dotsc,m\}$: then the two
relations
\[
\sum_{i=1}^s a_iw_i +\sum_{k\in I^r} \lambda_kx_k-\sum_{j\in J^r}
\lambda_jy_j =0, \qquad
\sum_{i=1}^s a_iw_i +\sum_{j\in J^r} \lambda_jx_j-\sum_{k\in I^r}
\lambda_ky_k =0
\]
correspond respectively to the numerical classes of the invariant
curves in $X$:
\[
V(\langle w_3,\ld,w_s,x_k,y_j \rangle_{k\in I^r,\,j\in J^r}), \qquad
V(\langle w_3,\ld,w_s,x_j,y_k \rangle_{j\in J^r,\,k\in I^r}).
\]
Since $X$ is Fano, they both have positive degree, so 
$\sum_ia_i\geq 1+|\sum_{k\in I^r}\lambda_k
-\sum_{j\in J^r}\lambda_j|\geq 1$.
\end{proof}

\noindent \emph{Step 6:} Suppose that $y_0=-x_0\in H'\cap G(\fx)$. 
Then every cone
$\sigma\in\fx$ decomposes as $\sigma=\eta+\tau$, with $\eta\subset
H'$, $\dim\eta=m$
and $\tau\cap H'=\{0\}$; $H'$ contains the fan of a Del
Pezzo variety $V^m$ and $X$ is a toric bundle over $W$ with fiber $V^m$. 
\begin{proof}[Proof of step 6]
Apply step 2 twice considering first $x_0$ and then $y_0$: 
it is immediate to see that $H'$ contains the fan of the Del Pezzo
variety $V^m$. Consider now a cone $\tau\in\fx$ such that $\tau\cap
H'=\{0\}$. By lemma~\ref{dopo}, we know that $\tau+\eta\in\fx$ for at
least one $m$-dimensional $\eta\subset H'$. Hence, again by step 2,
$\tau+\eta'\in\fx$ for all  $\eta'\subset H'$. Thus by
theorem~\ref{soleil} $X$ is a toric bundle over a smooth toric
variety  with fiber $V^m$. Since $W$ is an invariant section of the
bundle, it is isomorphic to the basis of the bundle.
\end{proof}

\noindent\emph{Step 7:}
suppose that $H'\cap
G(\fx)=\{x_0,x_1,y_1,\ld,x_m,y_m\}$.
Then there exist a smooth projective toric variety $X'$ 
with $G(\f_{X'})=G(\fx)$, such that $X'$ is a toric bundle
  over $W$ with fiber the pseudo Del Pezzo variety $\w{V}^m$.
\begin{proof}[Proof of step 7]
We remark that by step 4, all generators in
$G(\fx)\smallsetminus H'$ induce a  generator in $\f_{W}$.
The fan $\f_{X'}$ is defined as follows: 
\begin{enumerate}[1.]
\item $G(\f_{X'})=G(\fx)$;
\item the set
$\{\sigma\in\f_{X'}\,|\, \sigma\subset H'\}$ is the 
fan of the pseudo Del Pezzo variety $\w{V}^m$ in $H'$, with generators
$\{x_0,x_1,y_1,\ld,x_m,y_m\}$;
\item for every $\tau\in\fx$ such that $\tau\cap H'=\{0\}$ and
$\tau+\eta\in\fx$ for $\eta$ in \eqref{coni}, and for every
$\sigma\subset H'$ in the fan of $\w{V}^m$, we set $\tau+\sigma\in\f_{X'}$.
\end{enumerate}
By construction, $X'$ is a toric bundle with fiber $\w{V}^m$ (see
theorem~\ref{soleil}). The invariant subvarieties $V(\sigma)\subset
X'$ for $\sigma\in\f_{X'}(m)$, $\sigma\subset H'$, are invariant
sections of the bundle. If $\sigma$ is one of the cones in
\eqref{coni}, then $X\dasharrow X'$ is an isomorphism on
$V(\sigma)\simeq W$. Therefore the basis of the bundle is actually
$W$.  
\end{proof}

Let's consider now the linear subspace $H$ of $N_{\Q}$ spanned by all
pairs of symmetric generators in $\fx$. We want to prove that $H$
decomposes as a direct sum of subspaces as $H'$. 

Fix $h=\dim H$ pairs
$x_i+y_i=0$ that span $H$. If $G(\fx)\cap
H=\{x_1,y_1,\dotsc,x_{h},y_{h}\}$, 
then the desired decomposition is just $H=K$. 
If there is a generator $x_0\in H$ different from $\{x_i,y_i\}$, by
step 1, up to renaming $\{x_i,y_i\}$, we can suppose that
$x_0=y_1+\cdots+y_m$ with $2\leq m \leq h$.
Let $H'$ be the linear span of $\{x_1,\dotsc,x_m\}$. 
We have to show that if $v$ is another generator
in $H$ such that $v\neq -x_0$, then the smallest linear subspace $H''$
containing $v$ satisfies $H'\cap H''=\{0\}$.

Again by step 1, we know that 
$v=\sum_{i=1}^{h} \varepsilon_ix_i$ with
  $\varepsilon_i\in \{1,0,-1\}$ for all $i$.
Up to renaming $\{x_{m+1},y_{m+1},\dotsc,x_h,y_h\}$, we can assume that
\[ v= \sum_{i\in I}x_i + \sum_{j\in J}y_j+\sum_{k=m+1}^sy_k\]
where $I, J\subset\{1,\dotsc,m\}$ disjoint and $m<s\leq h$. Remark that it
must be $s>m$, because by step 3 we know that 
$v\not\in H'$ (and $x_0\not\in H''$).

Suppose first that $\#I\leq r$ and $\#J\leq r$ and 
choose a partition $I^r,J^r$ of $\{1,\dotsc,m\}$ such that $I^r\supseteq I$, 
$J^r\supseteq J$. Set 
$\eta=\langle x_i,y_j\rangle_{i\in I^r,\,j\in
  J^r}$. By step 4, we have $\langle v\rangle+\eta\in\fx$, and since
$\{v,x_i,y_j,x_{m+1},\ld,x_{s-1}\}_{i\in I^r,\,j\in J^r}$ are
linearly independent, lemma~\ref{dopo} implies that we can reorder
$\{x_{m+1},\ld,x_{s-1}\}$ in such a way that
\[ \eta+\langle v,x_{m+1},\ld,x_{t},y_{t+1},\ld,y_{s-1}\rangle\in\fx,
\text{ with } m\leq t\leq s-1,\]
and by step 2 the same is true for $\eta'=\langle x_j,y_i\rangle_{j\in
  J^r,\,i\in I^r}$.
Now, the relations
\[\!\!\!\sum_{k=t+1}^s y_k-(v+\sum_{i\in I}y_i+\sum_{j\in
  J}x_j+\sum_{h=m+1}^tx_h)=0,\quad
 v+\sum_{h=m+1}^tx_h+x_s-(\sum_{i\in I}x_i+\sum_{j\in
  J}y_j+\sum_{k=t+1}^{s-1} y_k)=0
\]
are effective by lemma~\ref{lemma}, and since $X$ is Fano, they have
positive degree: thus we get $\#I+\#J\leq s+m-2t-2$ and $\#I+\#J\leq 
2t+2-m-s$. So $\#I+\#J\leq 0$ and $\#I=\#J=0$: 
then $v=\sum_{k=m+1}^sy_k$ and $H'\cap H''=\{0\}$.

Suppose now that either $\#I$ or $\#J$ is bigger than $r$. In this
case, we get a contradiction; the proof
is very similar to the preceeding case. Suppose that $\#I>r$:
we can write
\[ v=x_1+\cdots+x_r+\sum_{i\in I'}x_i +\sum_{j\in J}y_j
+\sum_{k=m+1}^sy_k, \]
where $I=\{1,\ld,r\}\cup I'$. Reasonning as before we get the same
cones in $\fx$; now we consider the relations:
\begin{gather*}
 \sum_{i\in I'}x_i+\sum_{k=t+1}^s y_k-(v+y_1+\cdots+y_r+\sum_{j\in
  J}x_j+\sum_{h=m+1}^tx_h)=0,\\
 v+\sum_{i\in I'}y_i+\sum_{h=m+1}^tx_h+x_s-(x_1+\cdots+x_r+\sum_{j\in
  J}y_j+\sum_{k=t+1}^{s-1} y_k)=0.
\end{gather*}
Since they have positive degree, we get $\#I'\geq \#J+r$, hence
$\#J=0$ and $\#I'=r$. Thus $v=-x_0$, a contradiction. In the same way,
if $\#J>r$, we get $v=x_0$, again a contradiction.
\end{proof}

\begin{remark} In the same setting as step 7,
  it is not true in general that $H'$ contains the
fan of $\widetilde{V}^m $. Indeed, the sets $\{x_i,y_j\}_{i\in
  I^{r+s},\, j\in J^{r-s}}$ where $I^{r+s},J^{r-s}$ is a partition
of $\{1,\dotsc,m\}$ and $s\in\{1,\dotsc,r\}$, do not need to generate a cone
in $\fx$. 
\end{remark}

As an immediate corollary, we get a generalization of Ewald's
result~\cite{ewald}: 
\begin{prop}
\label{dipiu}
Let $X$ be a toric Fano variety and let $H$ be the linear subspace of
$N_{\Q}$ spanned by all pairs of symmetric generators in
$G(\fx)$. 

Suppose that $\dim H=n$.
 Then $X$ is a product
of Del Pezzo varieties, pseudo Del Pezzo varieties 
and a power of $\pr{1}$.

 Suppose that $\dim H=n-1$. Then there exists 
$X'\in\mathcal{B}_X$ which is a toric bundle over 
$\pr{1}$ with fiber a product
of Del Pezzo varieties, pseudo Del Pezzo varieties 
and a power of $\pr{1}$.
\end{prop}
In dimension 4, there are only two toric Fano varieties which have 3
linearly 
independent pairs of symmetric generators and such that
$X'$ is not isomorphic to $X$.
These 4-folds are $M1$ and $R3$, in Batyrev's notation
\cite{bat2}. In both cases $\ph\colon X'\dasharrow X$ is just a flip.

\begin{proof}[Proof of proposition \ref{dipiu}]
The statement for $\dim H=n$ is straightforward from theorem~\ref{cs}:
let $X\dasharrow X'\rightarrow
W$ be as in theorem~\ref{cs}, $(iii)$. Then 
$G(\f_W)$ consists of $\dim W$ linearly independent pairs of symmetric
generators
$\overline{x}_i,\overline{y}_i$, where $\overline{u}$ denotes the
 image of $u$ under the projection $N\rightarrow N/H'$, $H'=
 H_1\oplus\cdots\oplus H_a\oplus\widetilde{H}_1\oplus\cdots\oplus
 \widetilde{H}_b$.
 Since $W$ is Fano, we have
$W\simeq(\pr{1})^{\dim
  W}$. Moreover,  the primitive relations of $W$
$\overline{x}_i+\overline{y}_i=0$ remain unchanged when lifted to $N$,
hence $X'$ is actually the product $V^{h_1}\times\cdots\times
V^{h_a}\times\w{V}^{\w{h}_1}\times\cdots  
\times\w{V}^{\w{h}_b}\times(\pr{1})^{\dim W}$. Since $X'$ is Fano,
$X'\simeq X$. 

Assume $\dim H=n-1$ and fix $n-1$ pairs $x_i,y_i$ with
$i=1,\dotsc,n-1$ such that $x_1,\dotsc,x_{n-1}$ are a basis for $H$
and $x_i+y_i=0$ for all $i$.
\begin{claim} Let $v,w\in
G(\fx)\smallsetminus H$. Then either $\{v,w\}$ is a primitive
collection with relation $v+w-z=0$ and $z\in G(\fx)\cap H$, or
$\langle v,w\rangle\in\fx$ and 
$w=\sum_{i=1}^{n-1}\varepsilon_ix_i+\delta v$ with $|\delta|=1$ and
$\varepsilon_i\in\{1,0,-1\}$ for all $i$; moreover,
the number of non-zero $\varepsilon_i$ is odd.
\end{claim}
\begin{proof}[Proof of the claim]
If $v\in G(\fx)\smallsetminus H$, then by lemma~\ref{dopo}
$\{v,x_1,\dotsc,x_{n-1}\}$ is a basis of $N$. Hence, if $v,w\in
G(\fx)\smallsetminus H$, we have
$w=\sum_{i=1}^{n-1}\varepsilon_ix_i+\delta v$ with $|\delta|=1$.

Suppose that $\{v,w\}$ is a primitive collection. Since $v+w\neq 0$
(otherwise $\dim H=n$), the
associated primitive relation must be $v+w-z=0$, with $z\in G(\fx)$. 
Then $z=\sum_{i=1}^{n-1}\varepsilon_ix_i+(1+\delta)v$, hence it must
be $\delta=-1$ and $z\in H$. 

Let's assume that $\langle v,w\rangle\in\fx$. We can reorder
$\{x_1,\dotsc,x_{n-1}\}$ in such a way that
$w=\sum_{i=1}^r\varepsilon_ix_i+\delta v$ with $r\in\{1,\dotsc,n-1\}$
and $\varepsilon_i\neq 0$ for all $i=1,\dotsc,r$.
Fix $i\in\{1,\dotsc,r\}$. Since $\{v,w,x_1,\ld,\check{x}_i,\ld,x_r\}$
are linearly independendent, lemma \ref{dopo} implies that it is a
part of a basis of the lattice. Hence $|\varepsilon_i|=1$.

It remains to show that $r$ is odd. Up to renaming $\{x_i,y_i\}$, we
can assume that $w=\sum_{i=1}^r x_i+\delta v$, where
$|\delta|=1$. Again by lemma \ref{dopo}, there exists a partition
$I,J$ of $\{1,\dotsc,r-1\}$ such that $\langle
v,w,x_i,y_j\rangle_{i\in I,\,j\in J}\in\fx$. Consider now the two
relations:
\[ \sum_{i\in I}x_i +x_r+\delta v-w-\sum_{j\in J}y_j=0,\qquad
\sum_{j\in J}y_j+y_r+w-\delta v-\sum_{i\in I}x_i=0. \]
They are both effective by lemma \ref{lemma}, so they must have
positive degree: this gives $\#J=\#I+\delta-1$, hence $r=2\#I+\delta$
is odd. 
\end{proof}
We first prove the proposition in the case   $G(\fx)\cap
H=\{x_1,y_1,\dotsc,x_{n-1},y_{n-1}\}$ (namely, in the decomposition
given by theorem~\ref{cs}, we have $H=K$).
Then, we can restate the claim as follows: 
 for every $v,w\in
G(\fx)\smallsetminus H$, we have
$w=\sum_{i=1}^{n-1}\varepsilon_ix_i+\delta v$ with $|\delta|=1$ and
$\varepsilon_i\in\{1,0,-1\}$ for all $i$; moreover, 
the number of non-zero $\varepsilon_i$ is odd.

Remark that $\#(G(\fx)\smallsetminus H)\geq 2$, because since $X$ is
complete, there is at least one generator in both half-spaces cut by $H$. 
Let's show that $\#(G(\fx)\smallsetminus H)= 2$. By contradiction,
suppose that $G(\fx)\smallsetminus
H\supseteq\{v_1,v_2,w\}$ with $v_1,v_2$ in the same half-space cut by
$H$ and $w$ in the other one. By the claim
we have
$v_2=v_1+\sum_{i\in I}\varepsilon_ix_i$ and $w=-v_1
+\sum_{j\in J}\delta_jx_j$
with $I,J\subseteq\{1,\dotsc,n-1\}$, $|\varepsilon_i|=|\delta_j|=1$ for all
$i,j$, and $\#I$ and $\#J$ odd. Then
\[ w+v_2=\sum_{i\in
  (I\smallsetminus J)}\varepsilon_ix_i+\sum_{j\in(J\smallsetminus
  I)}\delta_j x_j+\sum_{i\in I\cap J}(\varepsilon_i+\delta_i)x_i. \]
Since by the claim
all the coefficients must be 1 or -1, we have
  $\varepsilon_i+\delta_i=0$ for all $i\in I\cap J$, and the non-zero
  coefficients are $\#(I\smallsetminus J)+\#(J\smallsetminus
  I)=\#I+\#J -2\#(I\cap J)$ an even number, a contradiction.

Therefore we have
$G(\fx)=\{x_1,y_1,\dotsc,x_{n-1},y_{n-1},v,w\}$ 
and up to renaming
$\{x_i,y_i\}$, we can assume $v+w=x_1+\cdots+x_r$,
$r\in\{1,\dotsc,n-1\}$ odd. Hence 
$G(\fx)=G(\f_{X'})$ where $X'$ has primitive relations:
$x_i+y_i=0$ for $i=1,\dots,n-1$, and
  $v+w-(x_1+\cdots+x_r)=0$.
Thus $X'$ is a $(\pr{1})^{n-1}$-bundle over $\pr{1}$. We can actually
say more: define $Y'$ to be the $(r+1)$-dimensional smooth toric variety
with $G(\f_{Y'})=\{x_1,y_1,\dotsc,x_r,y_r,v,w\}$ and primitive
relations
$x_i+y_i=0$ for $i=1,\dots,r$ and
  $v+w-(x_1+\cdots+x_r)=0$. 
Then $Y'$ is  a $(\pr{1})^{r}$-bundle over $\pr{1}$ and
$X'=(\pr{1})^{n-1-r}\times Y'$. Moreover, 
$X=(\pr{1})^{n-1-r}\times Y$, where
$Y$ is Fano and $G(\fy)=G(\f_{Y'})$.

We have now to prove the general case. Let $X\dasharrow X'\rightarrow
W$ be as in theorem~\ref{cs}, $(iii)$. Then $W$ is Fano and has at
least $\dim
W-1$ linearly independent pairs of symmetric generators
$\overline{x}_i+\overline{y}_i=0$, where $\overline{u}$ denotes the
image of $u$ under the projection $N\rightarrow N/H'$, $H'=
H_1\oplus\cdots\oplus H_a\oplus\widetilde{H}_1\oplus\cdots\oplus
\widetilde{H}_b$.
 By what preceeds, $G(\f_W)$ has exactly 
$2\dim W$ elements: the last two are $\overline{v}$, $\overline{w}$
such that 
$\overline{v}+\overline{w}=\overline{x}_1+\cdots+\overline{x}_r$ 
 for some $r\in\{0,\dotsc,\dim W\}$ (if $r>0$ then it must be odd).
In $X'$, by hypothesis we still have $x_i+y_i=0$; the last relation
lifts as $v+w=x_1+\cdots+x_r+u_1+\cdots+u_s$, where $\{u_1,\dotsc,u_s\}$
are contained in $H'$ (and again $r+s$ is odd). 
Hence $X'$ is a toric bundle over
$\pr{1}$ with fiber
 $V^{h_1}\times\cdots\times V^{h_a}\times\w{V}^{\w{h}_1}\times\cdots 
\times\w{V}^{\w{h}_b}\times(\pr{1})^{\dim W-1}$.
\end{proof}
\begin{remark}
We think that it should be possible, with a detailed analysis, to
classify all the possible $X$ Fano 
with $\dim H=n-1$, 
namely: to list all the possible decompositions of $v+w$ in $H$, or in
other words again, to list all the possible $X'\rightarrow \pr{1}$
with fiber a product
of Del Pezzo varieties, pseudo Del Pezzo varieties 
and a power of $\pr{1}$, such that there exists $X$ Fano with
$G(\fx)=G(\f_{X'})$. In the proof, we did this when the fiber of
$X'\rightarrow \pr{1}$ is just a power of $\pr{1}$. 
\end{remark}

\bigskip

\noindent{\bf Acknowledgements.} 
I would like to thank Laurent Bonavero and Olivier Debarre, for many
useful conversations about toric Fano varieties.

I am also very grateful to Lucia Caporaso, my advisor, 
for her constant guide.

\bigskip

\small

\begin{thebibliography}{BCDD02}

\bibitem[Bat82]{bat3}
Victor~V. Batyrev.
\newblock Toroidal {F}ano 3-folds.
\newblock {\em Mathematics of the USSR Izvestiya}, 19:13--25, 1982.

\bibitem[Bat99]{bat2}
Victor~V. Batyrev.
\newblock On the classification of toric {F}ano 4-folds.
\newblock {\em Journal of Mathematical Sciences (New York)}, 94:1021--1050,
  1999.

\bibitem[BCDD02]{mukai}
Laurent Bonavero, Cinzia Casagrande, Olivier Debarre, and St{\'e}phane Druel.
\newblock Sur une conjecture de {M}ukai.
\newblock Pr{\'e}publication 566, Institut Fourier, 2002.
\newblock {Preprint math.AG/0204314}.

\bibitem[BCW02]{bonwisncamp}
Laurent Bonavero, Fr{\'e}d{\'e}ric Campana, and Jaroslaw~A. Wi{\'s}niewski.
\newblock Vari{\'e}t{\'e}s projectives complexes dont l'{\'e}clat{\'e}e en un
  point est de {F}ano.
\newblock {\em Comptes Rendus de l'Acad{\'e}mie des Sciences}, 334:463--468,
  2002.

\bibitem[Bon00]{bonavero3}
Laurent Bonavero.
\newblock Toric varieties whose blow-up at a point is {F}ano.
\newblock {Preprint math.AG/0012229}, 2000.
\newblock To appear in Tohoku Mathematical Journal.

\bibitem[Cas01a]{contr}
Cinzia Casagrande.
\newblock Contractible classes in toric varieties.
\newblock {Preprint math.AG/0111332}, 2001.
\newblock To appear in Mathematische Zeitschrift.

\bibitem[Cas01b]{cras}
Cinzia Casagrande.
\newblock On the birational geometry of toric {F}ano 4-folds.
\newblock {\em Comptes Rendus de l'Acad{\'e}mie des Sciences}, 332:1093--1098,
  2001.

\bibitem[Cas01c]{fano}
Cinzia Casagrande.
\newblock Toric {F}ano varieties and birational morphisms.
\newblock {Preprint math.AG/0112007}, 2001.

\bibitem[Ewa88]{ewald}
G{\"u}nter Ewald.
\newblock On the classification of toric {F}ano varieties.
\newblock {\em Discrete Computational Geometry}, 3:49--54, 1988.

\bibitem[Ewa96]{ewald2}
G{\"u}nter Ewald.
\newblock {\em Combinatorial Convexity and Algebraic Geometry}, volume 168 of
  {\em Graduate Texts in Mathematics}.
\newblock Springer-Verlag, 1996.

\bibitem[Ful93]{fulton}
William Fulton.
\newblock {\em Introduction to Toric Varieties}.
\newblock Number 131 in Annals of Mathematics Studies. Princeton University
  Press, 1993.

\bibitem[Oda88]{oda}
Tadao Oda.
\newblock {\em Convex Bodies and Algebraic Geometry - An Introduction to the
  Theory of Toric Varieties}.
\newblock Springer-Verlag, 1988.

\bibitem[Rei83]{reid}
Miles Reid.
\newblock Decomposition of toric morphisms.
\newblock In {\em Arithmetic and Geometry, vol. II: Geometry}, number~36 in
  Progress in Mathematics, pages 395--418. Birkh{\"{a}}user, 1983.

\bibitem[Sat00]{sato}
Hiroshi Sato.
\newblock Toward the classification of higher-dimensional toric {F}ano
  varieties.
\newblock {\em Tohoku Mathematical Journal}, 52:383--413, 2000.

\bibitem[Sat02]{sato2}
Hiroshi Sato.
\newblock Toric {F}ano varieties with divisorial contractions to curves.
\newblock {Preprint math.AG/0205298}, 2002.

\bibitem[VK85]{VK}
V.~E. Voskresenski{\u{\i}} and Alexander Klyachko.
\newblock Toroidal {F}ano varieties and roots systems.
\newblock {\em Mathematics of the USSR Izvestiya}, 24:221--244, 1985.

\bibitem[Wi{\'s}02]{torimori}
Jaroslaw~A. Wi{\'s}niewski.
\newblock Toric {M}ori theory and {F}ano manifolds.
\newblock In {\em Geometry of Toric Varieties}, volume~6 of {\em S{\'e}minaires
  et Congr{\'e}s}, pages 249--272. Soci{\'e}t{\'e} Math{\'e}matique de France,
  2002.

\bibitem[WW82]{wat}
Keiichi Watanabe and Masayuki Watanabe.
\newblock The classification of {F}ano 3-folds with torus embeddings.
\newblock {\em Tokyo Journal of Mathematics}, 5:37--48, 1982.

\end{thebibliography}

\end{document}